\documentclass[12pt]{article}
\usepackage{amsmath,amssymb,latexsym,graphicx}
\newtheorem{theorem}{Theorem}
\newtheorem{lemma}[theorem]{Lemma}
\newtheorem{remark}[theorem]{Remark}
\newtheorem{corollary}[theorem]{Corollary}

\newtheorem{conjecture}[theorem]{Conjecture}
\newcommand{\RR}{\mathbb{R}}

\begin{document}
\author{Peter Kuchment\footnote{Mathematics Department,
Texas A\& M University, College Station, TX 77843-3368, USA.
kuchment@math.tamu.edu} and Leonid Kunyansky\footnote{Mathematics Department,
University of Arizona, Tucson, AZ 77843-3368, USA.
leonk@math.arizona.edu}}
\title{Mathematics of thermoacoustic tomography}
\maketitle
\begin{abstract}
The paper presents a survey of mathematical problems, techniques, and
challenges arising in the Thermoacoustic (also called Photoacoustic or
Optoacoustic) Tomography.
\end{abstract} 
\section{Introduction} 

Computerized tomography has had a huge impact on medical diagnostics.
Numerous
methods of tomographic medical imaging have been developed and are being
developed (e.g., the ``standard'' $X$-ray, single-photon emission,
positron emission, ultrasound, magnetic resonance, electrical impedance,
optical)
\cite{He1, Kak, biomed, Natt_old, Natt_new}.
The designers of these modalities strive to increase the image
resolution and contrast, and at the same time to reduce the costs and
negative health effects of these techniques. However, these goals are
usually rather contradictory. For instance, some cheap and safe methods
with good contrast (like optical or electrical impedance tomography)
suffer from low resolution, while some high resolution methods (such as
ultrasound imaging) often do not provide good contrast. Recently
researchers have been developing novel hybrid methods that combine different
physical types of signals, in hope to alleviate the deficiencies of
each of the types, while taking advantage of their strengths. The most
successful example of such a combination is the {\bf Thermoacoustic
Tomography (TAT)} \footnote{TAT is also called Photoacoustic (PAT)
or Optoacoustic (OAT) Tomography and is sometimes abbreviated as TCT,
which stands for Thermoacoustic Computed Tomography}\cite{Kruger}.
Albeit not being a common feature in clinics yet, TAT scanners are actively
researched, developed and already manufactured, for instance
by OptoSonics, Inc. (http://www.optosonics.com/), founded by the pioneer of
TAT R.~Kruger.

After a substantial effort, major breakthroughs have been achieved in the
last couple of years in the mathematical modeling of TAT. The aim of this
article is to survey this recent progress and to describe the relevant
models, mathematical problems, and reconstruction procedures arising in TAT,
and to provide references to numerous research publications on this topic.

The main thrust of this text is toward mathematical methods;
considerations of the text length, as well as authors' background do not
let us discuss in any detail industrial and physical set-ups and
parameters of the TAT technique, and limitations of the
corresponding mathematical models.
Fortunately, the excellent recent surveys by M.~Xu and L.-H.~V.~Wang
\cite{MXW_review} and by A.~A.~Oraevsky and A.~A.~Karabutov \cite{Oraev,
Oraev2} accomplish all of these tasks, and thus the reader is advised to
consult with them for all such issues (see also the recent textbook
\cite{Wang_book}). On the other hand, in spite of the significant
recent progress in mathematics of TAT, there is no comprehensive survey text
addressing in details the relevant mathematical issues, although the surveys
\cite{Oraev2,MXW_review} do mention some mathematical reconstruction
techniques.

The structure of the paper is a follows: Section \ref{S:TAT} contains a brief
description of the TAT procedure.
The next Section \ref{S:MathTAT} provides the mathematical formulation of
the TAT problem. In general, it is formulated as an inverse problem for the
wave equation. However, in the case of the constant sound speed, it can be
also described in terms of a spherical mean operator (a spherical analog of
the Radon transform). The section also contains the list of natural
questions to be addressed concerning this model. These issues are addressed
then one by one in the following sections. In particular, Section
\ref{S:uniqueness} discusses uniqueness of reconstruction, i.e. the question
of whether the data collected in TAT is sufficient for recovery of the
information of interest. Albeit, for all practical purposes this issue is
resolved in Corollary \ref{C:closed}, we provide an additional discussion of
unresolved uniqueness problems, which are probably of more academic
interest.  Section \ref{S:reconstruction} addresses inversion formulas and
algorithms. In Section \ref{S:partial} effects of having only partial data
are discussed. Section \ref{S:range} contains results concerning the so
called range conditions, i.e. the conditions that all ideal data must
satisfy. Section \ref{S:remarks} provides additional remarks and discussions
of the issues raised in the previous sections. The paper ends with an
Acknowledgments section and bibliography. Concerning the latter, we need to
mention that the engineering and biomedical literature on TAT is rather vast
and no attempt has been made in this text to create a comprehensive
bibliography of the topic from the engineering prospective. The references
in \cite{Oraev, Oraev2, BiomDiagn, BiomPhot, MXW_review} to a large extent fill this gap. The
authors, however, have tried to present a sufficiently complete review of
the existing literature on mathematics of TAT.


\section{Thermoacoustic tomography}\label{S:TAT}

\begin{figure}[!ht]
\begin{center}
\scalebox{0.7}{\includegraphics{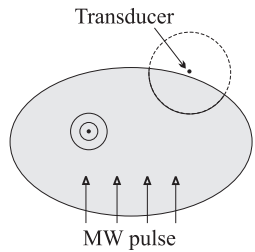}}
\end{center}
\caption{The TAT procedure.}
\label{F:tat}
\end{figure}

In TAT, a short duration EM pulse is sent
through a biological object (e.g., woman's breast in mammography) with the
aim of triggering a thermoacoustic response in the tissue.
As it is explained in \cite{MXW_review}, the radiofrequency (RF) and
the visible light frequency ranges are currently considered to be the
most suitable
for this purpose. Since mathematics works exactly the same way in both of
these frequency ranges, we will not make such distinction and will be
talking about just ``an EM pulse''. E.g., in Figure \ref{F:tat} a microwave
pulse is assumed.
\noindent
In most cases the pulse is spatially wide, so that the whole object is
more or less uniformly irradiated. Some part of EM energy is absorbed
throughout the object.
The amount of energy absorbed at a location strongly depends on
local biological properties of the cells. Oxygen saturation,
concentration of hemoglobin, density of the microvascular network
(angiogenesis), ionic conductivity, and water content are among the
parameters that influence the absorption strongly \cite{MXW_review}.
Thus, if the energy absorption distribution function $f(x)$ were known,
it would provide a great diagnostic tool. For instance, it could be useful
for detecting cancerous cells that absorb several times more energy
in the RF range than the healthy ones
\cite{Kruger, Oraev2, MXW1, MXW_review}.
However,  as an imaging tool neither RF waves, nor visual light alone
would provide acceptable resolution.  In the RF case, this
is due to the long wave length. One can use shorter microwaves, but
this will be at the expense of the penetration depth. In the optical
region, the problem is with the multiple scattering of light. So, a
different mechanism, the so called {\em Photoacoustic Effect} \cite{Gusev,
Tam, Wang_book, MXW_review}, is used
to image $f(x)$. Namely, the EM energy absorption results in
thermoelastic expansion and thus in a pressure wave $p(x,t)$ (an
ultrasound signal) that can be measured by transducers placed around
the object. Now one can attempt to recover the function $f(x)$ (the
image) from the measured data $p(x,t)$. Such a measuring scheme,
utilizing two types of waves, brings about the high resolution of
the ultrasound diagnostics and the high contrast of EM waves.
It overcomes the adverse effect of the low contrast of ultrasound
with respect to soft tissue. In fact, such a low contrast is a
good thing here, allowing one to assume in the first approximation that
the sound speed is constant. This often used approximation is not
always appropriate, but it is the most studied case at the moment.
Later on in this text we will describe some initial considerations of
the variable sound speed case, following \cite{Kuch_tat, JinWang}.

For this TAT method (and in particular, for the mathematical model
described below) to work, several conditions must be met. For instance,
the time duration of the EM pulse must be shorter than the time it
takes the sound wave to traverse the smallest feature that needs to be
reconstructed. The ultrasound detector must be able to resolve the time
scale of the duration of the EM pulse. On the other hand, the
transducer must be also able to detect much lower frequencies. Thus,
one needs to have extra-wide-band transducers, and these are currently
available. One can find the technical discussion of all these issues,
for instance, in \cite{Oraev2, MXW_review}. In this text we will assume
that all these conditions are met and thus the mathematical models
described are applicable.

In the next section we present a mathematical description of the
relation between $f(x)$ and $p(x,t)$ (similar mathematical problems arise in sonar \cite{LQ} and radar \cite{NC} imaging, as well as in geophysics \cite{Beylkin}).

\section{Mathematical model of TAT:\\ wave equation and the spherical
mean transform}\label{S:MathTAT}

\subsection{The wave equation model}\label{SS:wave equation}

We assume that the ultrasound speed at location $x$ is equal to $c(x)$.
Then, modulo some constant coefficients that we will assume all to be
equal to $1$, the pressure wave $p(x,t)$ satisfies the following
problem for the standard wave equation \cite{Diebold,Tam, MXW1}:
\begin{equation}\label{E:wave}
\begin{cases}
    p_{tt}=c^2(x)\Delta_x p, t\geq 0, x\in\RR^3\\
    p(x,0)=f(x),\\
    p_t(x,0)=0
    \end{cases}
\end{equation}
The goal is to find, using the data measured by transducers, the
initial value $f(x)$ at $t=0$ of the solution $p(x,t)$.

In order to formalize what data is in fact measured, one needs to
specify what kind of transducers is used, as well as the geometry of
the measurement. By the geometry of the measurement we mean the
distribution of locations of transducers used to collect the data.

We briefly describe here the commonly considered measurement procedure,
which uses point detectors. Line and planar detectors have also been
suggested (see Section \ref{SS:linear and planar}). It is too early to
judge which one of them will become most successful, but the one using
point transducers has been more thoroughly studied mathematically and
experimentally, and thus will be mostly addressed in this article. In
this case, the transducers are assumed to be point-like, i.e. of
sufficiently small dimension. A transducer at time $t$ measures the
average pressure over its surface at this time, which for the small
size of the transducer can be assumed to be just the value of $p(y,t)$
at the location $y$ of the transducer. Dimension count shows
immediately that in order to have enough data for reconstruction of the
function $f(x)$, one needs to collect data from the transducers'
locations $y$ running over a surface $S$ in $\RR^3$. Thus, the data at
the experimentalist's disposal is the function $g(y,t)$ that coincides
with the restriction of $p(x,t)$ to the set of points $y\in S$.

Taking into account that the measurements produce the values $g(y,t)$
of the pressure $p(x,t)$ of (\ref{E:wave}) on $S\times\RR^+$, the set
of equations (\ref{E:wave}) extends to become

\begin{equation}\label{E:wave_data}
\begin{cases}
    p_{tt}=c^2(x)\Delta_x p, t\geq 0, x\in\RR^3\\
    p(x,0)=f(x),\\
    p_t(x,0)=0\\
    p(y,t)=g(y,t), y\in S\times\RR^+
    \end{cases}
\end{equation}

\begin{figure}[ht]
\begin{center}
\scalebox{0.7}{\includegraphics{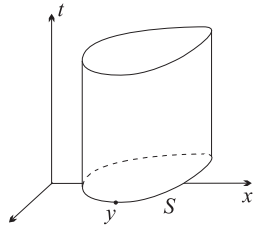}}
\label{F:cylinder}
\end{center}
\caption{An illustration to (\ref{E:wave_data}).}
\label{F:wave}
\end{figure}

The problem now reduces to finding the initial value $f(x)$ in
(\ref{E:wave_data}) from the knowledge of the lateral data $g(x,t)$
(see Figure \ref{F:cylinder}). A person familiar with PDEs might
suspect first that there is something wrong with this problem, since we
seem to have insufficient data for the recovery of the solution of the wave
equation in a cylinder from the lateral values alone. This, however, is
an illusion, since in fact there is a significant additional
restriction: the solution holds in the whole space, not just inside the
cylinder $S\times \RR^+$. We will see soon that in most cases, the data
is sufficient for recovery of $f(x)$.

\subsection{Spherical mean model}\label{SS:sphericalmean}

We now introduce an alternative formulation of the problem that works
in the constant speed case only. We will assume that the units are
chosen in such a way that $c(x)=1$. The known Poisson-Kirchhoff formula
\cite[Ch. VI, Section 13.2, Formula (15)]{CH} for the solution of
(\ref{E:wave}) gives

\begin{equation}\label{E:KP}
p(x,t)=c\frac{\partial}{\partial t}\left(t(Rf)(x,t)\right),
\end{equation}
where
\begin{equation}\label{E:mean}
(Rf)(x,r)= \int\limits_{|y|=1} f(x+ry)dA(y)
\end{equation}
is the  {\em spherical mean operator} applied to the function $f(x)$,
and $dA$ is the normalized area element on the unit sphere in $\RR^3$.
Hence, knowledge of the function $g(x,t)$ for $x\in S$ and all $t\geq
0$ essentially means knowledge of the spherical mean $Rf (x,t)$ at all
points $(x,t)\in S\times \RR^+$. One thus is lead to studying the
spherical mean operator $R:f\to Rf$ and in particular its restriction
$R_S$ to the points $x\in S$ only (these are the points where we place
transducers): \begin{equation}\label{E:Radon_S}
R_Sf(x,t)=\int\limits_{|y|=1}f(x+ty)dA(y), x\in S, t\geq 0.
\end{equation} This explains why, in many works on TAT, the spherical
mean operator has been the model of choice. Albeit the (unrestricted)
spherical mean operator has been studied rather intensively and for a
long time (e.g., \cite{Asgeirsson, CH, John}), its version $R_S$ with
the centers restricted to a subset $S$ appears to have been studied
since early 1990s only \cite{Agr}-\cite{AmbPatch}, \cite{And, Den,
Leon_Radon, Faw, FPR, Finch_even, FR, FR2, FR3, GGG, Gi, Kuc93, Kuch_AMS05,
KuchQuinto, Kunyansky,  LP1, LP2, LQ, Natt_new, Nils, Norton1, Norton2,
Palam_arc, Pal_book, Palam_funk, Patch, PopSush, PopSush2,
SchustQuinto, Zob} and offers quite a few new and often hard questions.

In what follows, we will alternate between these two (PDE and integral
geometry) interpretations of the TAT model, since each of them has its
own advantages.

\subsection{Main mathematical problems of TAT}\label{SS:math TAT}

We now formulate the typical list of problems one would like to address
in order to implement the TAT reconstruction.

\begin{enumerate}
\item For which sets $S\in \RR^3$ is the data collected by transducers
placed along $S$ sufficient for unique reconstruction of $f$? In terms
of the spherical mean operator, the question is whether $R_S$ has zero
kernel on an appropriate class of functions, say continuous with
compact supports.

\item If the data collected from $S$ is sufficient, what are inversion
formulas and algorithms?

\item How stable is the inversion?

\item What happens if the data is ``incomplete''?

\item What is the space of all possible ``ideal'' data $g(t,y)$
collected on a surface $S$? Mathematically (and in the constant sound
speed case) it is the question of describing {\em the range} of the
operator $R_S$ in appropriate function spaces. This question might seem
to be unusual (for instance, to people used to partial differential
equations), but in tomography importance of knowing the range of Radon
type transforms is well known. Such information is used to improve
inversion algorithms, complete incomplete data, discover and compensate
for certain data errors, etc. (e.g., \cite{Leon_Radon,GGG1, GGG,
GelfVil, Helg_Radon, Helg_groups, He1, Natt_old, Natt_new, Pal_book}).

\end{enumerate}

\section{Uniqueness of
reconstruction}\label{S:uniqueness}

Many of the problems of interest to TAT can be formulated in any
dimension $d$, albeit the practical dimensions are only $d=3$ and
$d=2$. We will consider an arbitrary dimension $d$ whenever we see this
suitable.

Let $S\subset \RR^d$ be the set of locations of transducers and $f$ be
a compactly supported function (one can show that for purposes of
uniqueness of reconstruction problem, one can always assume that $f$ is
smooth \cite{AQ}). Does the absence of the signal on the transducers,
i.e. $g(t,y)=0$ for all $t$ and $y$ in $S$, imply that $f=0$? If the
answer is a ``yes,'' we call $S$ - a {\bf uniqueness set}, otherwise a
{\bf non-uniqueness set}. In other words, in terms of TAT, the
uniqueness sets are those that distributing transducers along them
provides enough data for unique reconstruction of the function $f(x)$.

In terms of the wave equation, uniqueness sets are the sets of complete
observability, i.e. such that observing the motion on this set only,
one gets enough information to reconstruct the whole oscillation. In
terms of the spherical mean operator, the question is of whether the
equality $R_S f=0$ implies that $f=0$.

We will address this problem for the constant sound speed case first.

\subsection{Constant speed case}\label{SS:Unique_const}

As it has been discussed, the dimension count makes it clear that $S$
must be $(d-1)$-dimensional, i.e. a surface in $3D$ or a curve in $2D$.
We will also see that most of such surfaces are ``good'', i.e. are
uniqueness ones (or, in other words, provide enough information for
reconstruction). Thus, we should rather discuss the problem of
describing the ``bad'', non-uniqueness sets. The following simple
statement is very important and not immediately obvious.

\begin{lemma}\cite{AQ, LP1, LP2, Zob} Any non-uniqueness set $S$ is a
set of zeros of a (non-trivial) harmonic polynomial. In particular,
\begin{enumerate}
\item If there is no non-zero polynomial vanishing on $S$, then $S$ is
a uniqueness set.

\item If there is no non-zero harmonic function vanishing on $S$, then
$S$ is a uniqueness set.
\end{enumerate}
\end{lemma}

The proof of this lemma is very simple. It works under the assumption
of exponential decay of the function $f(x)$, not necessarily of
compactness of its support. It also introduces some polynomials that
play significant role in the whole analysis of the spherical mean
operator $R_S$.

Let $k\geq 0$ be an integer. Consider the convolution
\begin{equation}\label{E:Qk}
Q_k(x)=|x|^{2k}*f(x)=\int|x-y|^{2k}f(y)dy.
\end{equation}
This is clearly a polynomial of degree at most $2k$. Rewriting the
integral in polar coordinates centered at $x$ and using radiality of
$|x-y|$, one sees that $Q_k(x)$ is determined if we know the values
$Rf(x,t)$ of the spherical mean of $f$ centered at $x$:
$$
Q_k(x)=c_d\int\limits_0^\infty t^{2k+d-1}Rf(x,t)dt.
$$
In particular, If $R_S f\equiv 0$, then each polynomial $Q_k$ vanishes
on $S$.

Another observation that is easy to justify is that if the function $f$
is exponentially decaying (e.g., is compactly supported), then if all
polynomials $Q_k$ vanish identically, the function itself must be equal
to zero. (This is not necessarily true anymore if $f$ and its
derivatives decay only faster than any power, rather than
exponentially.)

Thus, we conclude that if $f$ is not identically equal to zero, then
there is at least one non-zero polynomial $Q_k$. Since, as we
discussed, equality $R_Sf=0$ implies that $Q_k|_S=0$, we conclude that
$S$ must be algebraic.

Now notice the following simple to verify equality (with a non-zero
constant $c_k$):
\begin{equation}\label{E:harm}
\Delta Q_k=c_kQ_{k-1},
\end{equation}
where $\Delta$ is the Laplace operator. This implies that the lowest
$k$ non-zero polynomial $Q_k$ is harmonic. Since $Q_k|_S =0$, this
proves the lemma.

Consider now the case when $S$ is a closed (hyper-)surface (i.e., the
boundary of a bounded domain). Since, as it is well known, there is no
non-zero harmonic function in the domain that would vanish at the
boundary (the spectrum of the Dirichlet Laplace operator is strictly
positive), we conclude that such $S$ is a uniqueness set for harmonic
polynomials. Thus, we get the following important
\begin{corollary}\label{C:closed}
\cite{AQ,Kuc93} Any closed surface is
uniqueness set for the spherical mean Radon transform.
\end{corollary}

An older alternative proof of this corollary provides an additional
insight into the problem. We thus sketch it here. Let us assume for
simplicity that the dimension $d\geq 3$ is odd (even dimensions require
a little bit more work). Suppose that the closed surface $S$ remains
stationary (nodal) for the oscillation described by (\ref{E:wave}).
Since the oscillation is unconstrained and the initial perturbation is
compactly supported, after a finite time, the interior of $S$ will
become stationary. On the other hand, we can think that $S$ is fixed
(since it is not moving anyway). Then, the energy inside $S$ must stay
constant. This is the contradiction that proves the statement of
Corollary \ref{C:closed}.

We will see in the next Section that the same method works in some
cases of variable sound speed, providing the needed uniqueness of
reconstruction result.

This corollary resolves the uniqueness problems for most practically
used geometries. It fails, however, if $f$ does not decay sufficiently
fast (see \cite{ABK}, where it is shown in which $L^p(\RR^d)$ classes
of functions $f(x)$ closed surfaces remain uniqueness sets).

It also provides uniqueness for some ``limited data'' problems. For
instance, if $S$ is an open (even tiny) piece of an analytic closed
surface $\Sigma$, it suffices. Indeed, if it did not, then it would be
a part of an algebraic non-uniqueness surface. Uniqueness of analytic
continuation would show then that the whole $\Sigma$ is a
non-uniqueness set, which we know to be incorrect. This result,
however, does not say that it would be practical to reconstruct using
observations from a tiny $S$. We will see later that this would not
lead to a satisfactory reconstructions, due to instabilities.

A geometry sometimes used is the planar one, i.e. detectors are placed
along a plane $S$ (line in the $2D$). In this case, there is no
uniqueness of reconstruction when the sound speed is constant. Indeed,
if $f(x)$ is odd with respect to $S$, then clearly all measured data
$g(t,y)$ will vanish. However, it is well known \cite{CH,John} that
functions even with respect to $S$ can be recovered. What saves the day
in TAT is that the object to be imaged is located on one side of $S$.
Then, extending $f(x)$ as an even function with respect to $S$, one can
still recover it from the data.

Although, for all practical purposes the uniqueness of reconstruction
problem is essentially resolved by the Corollary \ref{C:closed}, the
complete understanding of uniqueness problem has not been achieved yet.
Thus, we include below some known theoretical results and open
problems.

\subsubsection{Non-uniqueness sets in
$\mathbb{R}^2$.}\label{SS:2d_nonunique}

\begin{figure}[ht]
\begin{center}
\scalebox{0.7}{\includegraphics{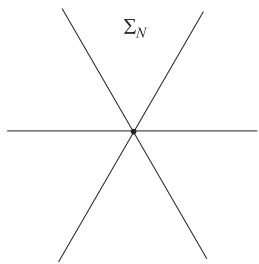}}
\end{center}
\caption{Coxeter cross $\Sigma_N$.}
\label{F:coxeter}
\end{figure}

In this Section, we follow the results and exposition of
\cite{AQ,LP1,LP2} in discussing uniqueness sets in $2D$. What are
simple examples of non-uniqueness sets? As we have already mentioned,
any line $S$ (or a hyperplane in higher dimensions) is a non-uniqueness
set, since any function $f$ odd with respect to $S$ will clearly
produce no signal: $R_Sf=0$. Analogously, consider a {\em Coxeter
system} $\Sigma _N$ of $N$ lines passing through a point and forming
equal angles (see Fig. \ref{F:coxeter}).

Choosing the intersection point as the pole and expanding functions
into Fourier series with respect to the polar angle, it is easy to
discover existence of an infinite dimensional space of functions that
are odd with respect to each of the $N$ lines. Thus, such a cross
$\Sigma _N$ is also a non-uniqueness set. Less obviously, one can use
the infinite dimensional freedom just mentioned to add any finite set
$\Phi$ of points still preserving non-uniqueness.
The following major and very non-trivial result was  conjectured in
\cite{LP1,LP2} and proven in \cite{AQ}. It shows that there are no
other bad sets $S$ besides the ones we have just discovered:

\begin{theorem}\label{T:AQ}
A set $S\subset \mathbb{R}^2$ is a non-uniqueness set for the spherical
mean transform in the space of compactly supported functions, if and
only if
$$S\subset \omega
\Sigma _N\cup \Phi,
$$
where $\Sigma _N$ is a Coxeter system
of lines, $\omega $ is a rigid motion of the plane, and $\Phi$ is a
finite set.
\end{theorem}

A sketch of a rather intricate proof of this result is provided in
Section \ref{SS:remunique}.

\subsubsection{Higher dimensions}\label{SS:higherdim}

Here we present a believable conjecture of how the result should look
like in higher dimensions.

\begin{conjecture}\label{C:n-dim}\cite{AQ}A set $S\subset \mathbb{R}^d$ is
a non-uniqueness set if and only if $ S\subset \omega \Sigma \cup \Phi
, $ where $\Sigma$ is the surface of zeros of a homogeneous harmonic
polynomial, $\omega $ is a rigid motion of $\RR^d$, and $\Phi$ is an
algebraic surface of dimension at most $d-2$.
\end{conjecture}

\begin{figure}[!ht]
\begin{center}
\scalebox{0.5}{\includegraphics{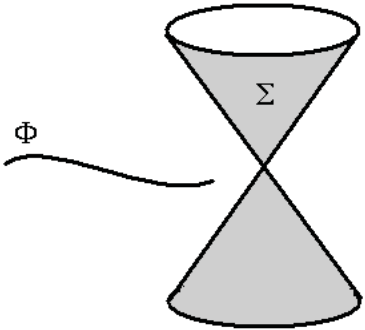}}
\end{center}
\caption{A picture of a $3$-dimensional non-uniqueness set.}
\label{F:n-dim}
\end{figure}
The progress towards proving this conjecture has been slow, albeit some
partial cases have been treated (\cite{Agr}-\cite{AmbKuc_inj}). E.g.,
in some cases one can prove that $S$ is a ruled surface (i.e., consists
of lines), but proving that these lines (rules) pass through a common point
remains a challenge. It is known, though, that both the zero sets of
homogeneous harmonic polynomials and algebraic subsets of dimension at
most $d-2$ are non-uniqueness sets \cite{Agr2, AQ}, and thus one should
avoid using them as placements of transducers for TAT.

\subsubsection{Relations to other areas of analysis}\label{SS:relations}

The problem of injectivity of $R_S$ has relations to a wide
variety of areas of analysis (see \cite{Agr, AQ} for many examples). In
particular, the following interpretation is important:

\begin{theorem} \cite{AQ,Kuc93} The following statements are equivalent:
\begin{enumerate}
\item $S \subset \mathbb{R}^d$ is a non-uniqueness set for the spherical
mean operator.\\
\item $S$ is a nodal set for the wave equation, i.e. there exists a
non-zero compactly supported $f$ such that the solution of the wave
propagation problem
$$
\begin{cases}
\frac{\partial ^2u}{\partial t^2}=\Delta u,\\
u(x,0)=0,\\
u_t(x,0)=f(x) \\
\end{cases}
$$
vanishes on $S$ for any moment of time.\\
\item $S$ is a nodal set for the heat equation, i.e. there exists a
non-zero compactly supported $f$ such that the solution of the problem
$$
\begin{cases}
\frac{\partial u}{\partial t}=\Delta u,\\
u(x,0)=f(x)
\end{cases}
$$
vanishes on $S$ for any moment of time.
\end{enumerate}
\end{theorem}
The interpretation in terms of the wave equation provides important PDE
tools and insights, which have lead to a recent progress
\cite{FPR,AmbKuc_inj} (albeit it has not lead yet to a complete
alternative proof of Theorem \ref{T:AQ}). The rough idea, originally
introduced in \cite{FPR}, is that if $S$ is a nodal set, then it might
be considered as the fixed boundary. In this case, the signals must go
around $S$. However, in fact, there is no obstacle, so signals can
propagate along straight lines. Thus, in order to avoid discrepancies
in arrival times, $S$ must be very special. One can find details in
\cite{FPR} and in \cite{AmbKuc_inj}.

\subsection{Uniqueness in the case of a variable sound
speed}\label{SS:uniqvariable}

It is shown in \cite[Theorem 4]{FR2} that uniqueness of reconstruction
also holds in the case of a smoothly varying (strictly positive) sound
speed, if the source function $f(x)$ is completely surrounded by the
observation surface $S$ (in other words, if there is no US signal
coming from outside of $S$). The proof uses the celebrated unique
continuation result by D.~Tataru \cite{Tataru}.

One can also establish uniqueness of reconstruction in the case of the
source not necessarily completely surrounded by $S$. However, here we
need to impose an additional non-trapping condition on the sound speed.
We assume that the sound speed is strictly positive $c(x)>c>0$ and such
that $c(x)-1$ has compact support, i.e. $c(x)=1$ for large $x$.

Consider the Hamiltonian system in $\RR^{2n}_{x,\xi}$ with
the Hamiltonian $H=\frac{c^2(x)}{2}|\xi|^2$:
\begin{equation}
\begin{cases}
x^\prime_t=\frac{\partial H}{\partial \xi}=c^2(x)\xi \\
\xi^\prime_t=-\frac{\partial H}{\partial x}=-\frac 12 \nabla
\left(c^2(x)\right)|\xi|^2 \\
x|_{t=0}=x_0, \xi|_{t=0}=\xi_0.
\end{cases}
\end{equation}
The solutions of this system are called {\em bicharacteristics}
and their projections into $\RR^n_x$ are {\em rays}.

We will assume that the following {\bf non-trapping condition}
holds:\\
all rays (with $\xi_0\neq 0$) tend to infinity when $t \to
\infty$.

\begin{theorem}\cite{Kuch_tat}\label{T:uniqueness}
Under the non-trapping conditions formulated above, compactly
supported function $f(x)$ is uniquely determined by the data $g$
measured on $S$ for all times. (No assumption of $f$ being supported
inside $S$ is imposed.)
\end{theorem}

One should mention that ray trapping can occur for some sound speed
profiles. For instance, if $c(x)=|x|$ for some range $r_1<|x|<r_2$,
then there are rays trapped in this spherical shell. We are not sure
what happens in this case to the uniqueness of reconstruction statement
of Theorem \ref{T:uniqueness} and inversion formula of Theorem
\ref{T:reconstr_speed}.

\section{Reconstruction: formulas and examples}\label{S:reconstruction}

Here we will address the procedures of actual reconstruction of the
source $f(x)$ from the data $g(t,y)$ measured by transducers.

\subsection{Constant sound speed}\label{SS:const_reconstr}

We assume here that the sound speed is constant and normalized to be
equal to $1$.

\subsubsection{Inversion formulas}\label{SS:inversion}

Before we move to our case of interest, which is spheres centered on a
closed surface $S$ surrounding the object to be imaged, we briefly refer to
related but somewhat different work. Namely, the problem of recovering
functions from integrals over spheres centered on a (hyper)plane $S$
has attracted a lot of attention over the years. Albeit, as it has been
mentioned before, there is no uniqueness in this case (functions odd
with respect to $S$ are annihilated), even functions can be recovered.
Thus also functions supported on one side of the plane can be as well,
by means of their even extension. Many explicit inversion formulas and
procedures have been obtained for this situation \cite{And, Den, Faw,
GGG, Gi, Kostli, Natt_new, Nils, Palam_arc, Pal_book, Ramm, Ramm2, Rom}. We will not
provide any details here, since this acquisition geometry is not very
useful. In particular, this is due to ``invisibility'' of some parts of
the interfaces, see Section \ref{S:partial}, which arises from
truncating the plane. The same problem is encountered with some other
unbounded acquisition surfaces, such as a surface of an ``infinitely''
long cylinder.

Thus, it is more practical to place transducers along a closed surface
surrounding the object. The simplest surface of this type is a sphere.

\subsubsection{Fourier expansion methods}\label{SS:Fourier}

Let us assume that $S$ is the unit sphere in $\RR^n$.
We would like to reconstruct a function $f(x)$ supported inside $S$
from the known values of its spherical integrals $g(y,r)$ with the
centers on $S$:
\begin{equation}
g(y,r)=\int\limits_{\mathbb{S}^{n-1}}f(y+r\mathbf{\omega}
)r^{n-1}d\mathbf{\omega},\qquad y\in S.  \notag
\end{equation}
The first inversion procedures for the case of spherical acquisition
were described in \cite{Norton1} in $2D$ and in \cite{Norton2} in $3D$.
These solutions were obtained by harmonic decomposition of the measured
data and the sought function, and by equating coefficients of the
corresponding Fourier series.

In particular, the 2-D algorithm of~\cite{Norton1} is based on the
Fourier decomposition of $f$ and $g$ in angular variables:
\begin{equation}
f(x)=\sum_{-\infty}^{\infty}f_{k}(\rho)e^{ik\varphi},
\quad x = (\rho \cos (\varphi), \rho \sin (\varphi) )
\label{nortonfourf}
\end{equation}
\[
g(y(\theta),r)=\sum_{-\infty}^{\infty}g_{m}(r)e^{ik\theta},
\quad y = (R \cos (\theta), R \sin (\theta) ).
\]
Following~\cite{Norton1} we consider the Hankel transform
$\hat{g}_{m,J}(\lambda)$ of the Fourier coefficients $g_{m}(r)$ (divided by
$2\pi r$)
\begin{equation}
\hat{g}_{m,J}(\lambda)=\int_{0}^{2R}g_{m}(r)J_{0}(\lambda r)dr=\mathcal{H}
_{0}\left(  \frac{g_{m}(r)}{2\pi r}\right)  .\label{besselintegr}
\end{equation}
To simplify the presentation we introduce the convolution $G_{J}(\lambda,y)$
of the sought function with the Bessel function $J_{0}(\lambda|x-y|).$
\begin{equation}
G_{J}(\lambda,y)=\int\limits_{\Omega}f(x)J_{0}(\lambda
|x-y|)dx,\label{nortongj}
\end{equation}
One can notice that $\hat{g}_{m,J}(\lambda)$ are the Fourier coefficients of
$G_{J}(\lambda,y)$ in $\theta$:
\begin{equation}
\hat{g}_{m,J}(\lambda)=\frac{1}{2\pi}\int_{0}^{2\pi}G_{J}(\lambda
,y)e^{-im\theta}d\theta.\label{nortonfourg}
\end{equation}
Now coefficients $f_{m}(\rho)$ can be recovered from $g_{m}(r)$ by application
of the addition theorem for the Bessel function $J_{0}(\lambda|x-y|)$:
\begin{equation}
J_{0}(\lambda|x-y|)=\sum_{-\infty}^{\infty}J_{m}(\lambda|x|)J_{m}
(\lambda|y|)e^{-im(\varphi-\theta)}.\label{nortonaddj}
\end{equation}
Indeed, by substituting equations (\ref{nortonfourf}) and (\ref{nortonaddj})
into (\ref{nortongj}), and (\ref{nortongj}) into (\ref{nortonfourg}) one
obtains
\[
\hat{g}_{m,J}(\lambda)=2\pi J_{m}(\lambda|R|)\int_{0}^{2R}f_{m}(\rho
)J_{m}(\lambda\rho)\rho d\rho=\mathcal{H}_{m}(f_{m}(\rho)),
\]
where $\mathcal{H}_{m}$ is the $m$-th order Hankel transform. Since the latter
transform is self-invertible, the coefficients $f_{m}(\rho)$ can be recovered
by the following formula
\begin{equation}
f_{m}(\rho)=\mathcal{H}_{m}\left[  \frac{\hat{g}_{m,J}(\lambda)}{J_{m}
(\lambda|R|)}\right]  =\mathcal{H}_{m}\left(  \frac{1}{J_{m}(\lambda
|R|)}\mathcal{H}_{0}\left[  \frac{g_{m}(r)}{2\pi r}\right]  \right)
,\label{norton2dformula}
\end{equation}
which is the main result of~\cite{Norton1}. Function $f(x)$ can now be
reconstructed by summing series (\ref{nortonfourf}).

Note that the above method requires a division of the Hankel transform of
the measured data by Bessel functions $J_{m}$ that have infinitely many zeros.
Theoretically, there is no problem; the Hankel transform $\mathcal{H}
_{0}\left[  \frac{g_{m}(r)}{2\pi r}\right]  $ has to have zeros that would
cancel those in the denominator. However, since the measured data always
contain some error, the exact cancelation is not likely to happen, and one
needs a sophisticated regularization scheme to keep the total error bounded.

This problem can be avoided by replacing in (\ref{besselintegr}) Bessel
function $J_{0}$ by Hankel function $H_{0}^{(1)}$:

\[
\hat{g}_{m,H}(\lambda)=\int_{0}^{2R}g_{m}(r)H_{0}^{(1)}(\lambda r)dr.
\]
The addition theorem for $H_{0}^{(1)}(\lambda|x-y|)$ takes form
\[
H_{0}^{(1)}(\lambda|x-y|)=\sum_{-\infty}^{\infty}J_{m}(\lambda|x|)H_{m}
^{(1)}(\lambda|y|)e^{-im(\varphi-\theta)},
\]
and by proceeding as before one can obtain the following formula for $f_{m}
(\rho)$:
\[
f_{m}(\rho)=\mathcal{H}_{m}\left[  \frac{\hat{g}_{m,H}(\lambda)}{H_{m}
^{(1)}(\lambda|R|)}\right]  =\mathcal{H}_{m}\left(  \frac{1}{H_{m}
^{(1)}(\lambda|R|)}\int_{0}^{2R}g_{m}(r)H_{0}^{(1)}(\lambda r)dr\right)  .
\]
Unlike $J_{m}$, Hankel functions $H_{m}^{(1)}(t)$ do not have zeros for all
real values of $t$ and therefore problems with division by zeros do not arise
in this amended version of the method~\cite{Norton1}.

This derivation can be repeated in 3-D, with the exponentials $e^{ik\theta}$
replaced by the spherical harmonics, and with cylindrical Bessel functions
replaced by their spherical counterparts. By doing this one will arrive at
the Fourier series method of~\cite{Norton2}. Our use of Hankel function
$H_{0}^{(1)}$ above is similar to the way the authors of~\cite{Norton2}
utilized spherical Hankel function $h_{0}^{(1)}$ to avoid the divisions by
zero.

\subsubsection{Filtered backprojection methods}\label{SS:backproj}

The favorite way of inverting Radon transform for tomography purposes
is by using filtered backprojection type formulas, which involve
filtration in Fourier domain followed (or preceded) by a
backprojection. In the case of the set of spheres centered on a closed
surface (e.g., sphere) $S$, one expects such a formula to involve a
filtration with respect to the radius variable and then some
integration over the set of spheres passing through the point of
interest. For quite a while, no such type formula had been discovered.
This did not prevent practitioners from reconstructions, since good
approximate inversion formulas (parametrices) could be developed,
followed by an iterative improvement of the reconstruction, see e.g.
reconstruction procedures in \cite{PAT, MXW1, XFW,
XXW, XWAK}, and especially \cite{PopSush, PopSush2}.

The first set of exact inversion formulas of the filtered
backprojection type was discovered in \cite{FPR}. These formulas were
obtained only in odd dimensions. Several different variations of such
formulas (different in terms of exact order of the filtration and
backprojection steps) were developed. Let us denote by $g(p,r)=r^2R_S f$ the
spherical integral, rather than the average, of $f$. Then various versions
of the $3D$ inversion formulas that reconstruct a function $f(x)$ supported
inside $S$ from its the spherical mean data $R_S f$, read:
\begin{equation}\label{E:FPR3d}
\begin{array}{c}
f(x)=-\frac{1}{8\pi ^{2}R}\Delta _{x}\int\limits_{\partial
B}g(y,|y-x|)dA(y),\\
f(x)=-\frac{1}{8\pi ^{2}R}\int\limits_{\partial B}\left( \frac{d^{2}
}{dt^{2}}g(y,t)\right) \left. {\phantom{\rule{1pt}{8mm}}}\right|
_{t=|y-x|}dA(y),\\
f(x)=-\frac{1}{8\pi ^{2}R}\int\limits_{\partial B}\left( \frac{d}{dt
}\left( \frac{1}{t}\frac{d}{dt}\frac{g(y,t)}{t}\right) \right)
\left. {\phantom{\rule{1pt}{8mm}}}\right| _{t=|y-x|}dA(
y).
\end{array}
\end{equation}
Recently, analogous formulas were obtained for even dimensions in
\cite{Finch_even}. Denoting by $g$, as before the spherical integrals
(rather than averages) of $f$, the formulas of \cite{Finch_even} in
$2D$ look as follows:
\begin{equation}\label{E:Finch2D}
f(x)=\frac{1}{2\pi R}\Delta \int\limits_{\partial
B}\int\limits_{0}^{2R}g(y,t)\log (t^{2}-|x-y
|^{2})\ dt\ dl(y),
\end{equation}
or
\begin{equation}
f(x)=\frac{1}{2\pi R}\int\limits_{\partial
B}\int\limits_{0}^{2R}\frac{\partial }
{\partial t}\left( t\frac{\partial }
{\partial t}\frac{g(y,t)}{t}\right) \log (t^{2}-|x-y|^{2})\ dt\ dl(y),
\end{equation}
A different set of explicit inversion formulas that work in arbitrary
dimensions was presented in \cite{Kunyansky}.
\begin{equation}\label{E:kunyansky}
f(x)=\frac{1}{4(2\pi )^{n-1}}\mathrm{div}\int\limits_{\partial B}
\mathbf{n}(y)h(y,|x-y|)dA(y).
\end{equation}
Here
\begin{eqnarray}
h(y,t)& =\int\limits_{\mathbb{R}^{+}}\left[ Y(\lambda t)\left(
\int\limits_{0}^{2R}J(\lambda t^{\prime })g(y,t^{\prime
})dt^{\prime }\right) \right. \notag \\
& -J(\lambda t)\left. \left( \int\limits_{0}^{2R}Y(\lambda
t^{\prime})g( y,t^{\prime })dt^{\prime }\right) \right]
\lambda ^{2n-3}d\lambda ,       \label{genfilt}
\end{eqnarray}
\begin{equation*}
J(t) =\frac{J_{n/2-1}(t)}{t^{n/2-1}}, \quad \quad
Y(t) =\frac{Y_{n/2-1}(t)}{t^{n/2-1}},
\end{equation*}
$J_{n/2-1}(t)$ and $Y_{n/2-1}(t)$ are respectively the Bessel and
Neumann functions of order $n/2-1$, and $\mathbf{n}(y)$ is the vector
of exterior normal to $\partial B$.

In 2-D equations (\ref{E:kunyansky}), (\ref{genfilt}) can be simplified
to yield the following reconstruction formula
\begin{equation}
f(x)=-\frac{1}{2\pi^{2}}\mathrm{div}\int\limits_{\partial B}\mathbf{n}
(y)\left[  \int\limits_{0}^{2R}g(y,t^{\prime})\frac{1}{|x-y|^{2}-
t^{\prime2}
}dt^{\prime}\right]  dl(y).  \notag 
\end{equation}
A similar simplification is also possible in $3D$ resulting in the
formula
\begin{equation}
f(x)=\frac{1}{8\pi^{2}}\mathrm{div}\int\limits_{\partial B}\mathbf{n}
(y)\left(  \frac{1}{t} \frac{d}{dt}\frac{g(y,t)}{t}\right)  \left.
{\phantom {\rule
{1pt}{8mm}}}\right|  _{t=|y-x|}dA(y).        \label{E:Kunya3D}
\end{equation}
Equation (\ref{E:Kunya3D}) is equivalent to one of the formulas derived in
\cite{MXW2} for the $3D$ case. It is interesting to notice that the
``universal'' formula of \cite{MXW2} holds for all geometries when the
backprojection type formulas are known: spherical, cylindrical, and
planar. It is not very likely that such explicit formulas would be
available for any closed surfaces $S$ different from spheres (see a related
discussion in \cite{Anastasio, Beylkin}).

\subsubsection{Series solutions for arbitrary
geometries}\label{SS:series}

Although, as we have just mentioned, we do not expect such explicit
formulas to be derived for non-spherical closed surfaces $S$, there is,
however, a different approach \cite{Kun_series} that theoretically
works for any closed $S$ and that is practically useful in some
non-spherical geometries.

Let $\lambda _{m}^{2}$ and $u_{m}(x)$ be the eigenvalues and
normalized eigenfunctions of the Dirichlet Laplacian $-\Delta $ on the
interior $\Omega$ of a closed surface $S$:
\begin{align}
\Delta u_{m}(x)+\lambda _{m}^{2}u_{m}(x)& =0,\qquad
x\in \Omega ,\quad \Omega \subseteq \mathbb{R}^{n},  \label{Helmeq}
\\
u_{m}(x)& =0,\qquad x\in S ,  \notag \\
||u_{m}||_{2}^{2}& \equiv \int\limits_{\Omega }|u_{m}(x)|^{2}d
x=1.  \notag
\end{align}
As before, we would like to reconstruct a compactly supported function
$f(x)$ from the known values of its spherical integrals $g(y,r)$ with
the centers on $S$:
\begin{equation*}
g(y,r)=\int\limits_{\omega^{n-1}}f(y+r\mathbf{\omega}
)r^{n-1}d\mathbf{\omega},\qquad y\in S.  \notag
\end{equation*}
We notice that $u_{m}(x)$ is the solution of the Dirichlet problem
for the Helmholtz equation with zero boundary conditions and the wave
number $\lambda _{m}$, and thus it admits the
Helmholtz representation\begin{equation}
u_{m}(x)=\int_{\partial \Omega }\Phi _{\lambda _{m}}(|x-
y|)\frac{\partial }{\partial \mathbf{n}}u_{m}(y)ds(y)\qquad x\in \Omega
,  \label{helmdiscr}
\end{equation}
where $\Phi _{\lambda _{m}}(|x-y|)$ is a free-space
rotationally invariant Green's function of the Helmholtz
equation~(\ref{Helmeq}).

The eigenfunctions $\left\{
u_{m}(x)\right\} _{0}^{\infty }$ form an orthonormal basis in
$L_{2}(\Omega ).$
Therefore, $f(x)$ can be represented by the series
\begin{equation}
f(x)=\sum_{m=0}^{\infty }\alpha _{m}u_{m}(x)     \label{fourierser}
\end{equation}
with
\begin{equation}
\alpha _{m}=\int_{\Omega }u_{m}(x)f(x)dx.         \notag
\end{equation}
Since $f(x)$ is $C_{0}^{1},$ series~(\ref{fourierser}) converges
pointwise. A reconstruction formula of $\alpha_m$, and thus of $f(x)$,
will result if we substitute representation~(\ref{helmdiscr}) into
(\ref{fourierser}) and interchange the order of integrations. Indeed,
after a brief calculation we will get
\begin{equation}
\alpha _{m} =\int_{\Omega }u_{m}(x)f(x)dx  =\int_{\partial \Omega
}I(y,\lambda _{m})\frac{\partial }{ \partial \mathbf{n}}u_{m}(y)dA(x),
\label{serkoef1a}
\end{equation}
where
\begin{equation}
I(y,\lambda)=\int_{\Omega }\Phi _{\lambda}(|x-
y|)f(x)dx.  \label{Greenint}
\end{equation}

Certainly, the need to know the spectrum and eigenfunctions of the
Dirichlet Laplacian imposes a severe constraint on the surface $S$.
However, there are simple cases when the eigenfunctions are well known,
and fast summation formulas for the corresponding
series are available.  Such is the case of a cubic measuring
surface $S$ (see \cite{Kun_series}); the eigenfunctions $u_m$ are
products of sine functions
\begin{equation}
u_{m}(x) = \frac{8}{R^{3}}\sin \frac{\pi m_{1}x_{1}}{R}
\sin \frac{\pi m_{2}x_{2}}{R}\sin \frac{\pi m_{3}x_{3}}{R},
\end{equation}
where $m=(m_{1},m_{2},m_{3}),$ $m_{1},m_{2},m_{3}\in \mathbb{N}$,
and the eigenvalues are easily found as well
\begin{equation}
\lambda_m = \pi^{2} |m|^{2}/{R^2}.
\end{equation}
Sum (\ref{fourierser}) is just a regular 3-D Fourier sine series
easily computable by application of the Fast Sine Fourier transform
algorithm. The algorithmic trick that allows one to calculate fast the
coefficients $\alpha _m$ consists in computing first integrals
(\ref{Greenint}) on a uniform mesh in $\lambda$. This is easily done by
a one-dimensional Fast Cosine Fourier transform algorithm, with
$\Phi_\lambda(t) = \cos (\lambda t)/t$. The normal derivatives of
$u_{m}(x)$ are also products of sine functions, this time
two-dimensional ones. This, in turn, permits rapid evaluation of
integrals $\int_{\partial \Omega_i} I(y,\lambda)\frac{\partial }
{ \partial \mathbf{n}}u_{m}(y)dA(x)$ for each mesh value of $\lambda$,
and for each one of the six faces $\partial \Omega _{i}, i=1,...,6$
of the cube. Finally, the computation of $\alpha _{m}$ using
equation (\ref{serkoef1a}) reduces to the interpolation in the spectral
parameter $\lambda$, since the integrals in the right hand side of this
equation have been computed for the mesh values of this parameter
(not for $\lambda_m$). Due to oscillatory nature of the integrals
(\ref{Greenint}) a low order interpolation here would lead to inaccurate
reconstructions. Luckily, however, these integrals are analytic
functions of parameter $\lambda$ (due to the finite support of $g$). Hence,
high order polynomial interpolation is applicable, and numerics yields very
good results.

The algorithm we just described requires $\mathcal{O}(m^3 \log m )$
floating point operations if the reconstruction is to be performed on
an $m \times m \times m$ Cartesian grid, from comparably discretized
data measured on a cubic surface. In practical terms, it yields
reconstructions in the matter of several seconds on grids with total
number of nods exceeding a million~\cite{Kun_series}.

\subsubsection{Time reversal (backpropagation)
methods}\label{SS:backprop}

In the constant speed case, the following approach is possible in $3D$:
due to the validity of the Huygens' principle (i.e., the signal escapes
from any bounded domain in finite time), the pressure $p(t,x)$ inside
$S$ will become equal to zero for any time $T$ larger than the time
required to cross the domain (i.e., time that it takes the sound to
move along the diameter of $S$, which for $c=1$ equals the diameter).
Thus, one can impose the zero conditions on $p(t,x)$ for $t=T$ and
solve the wave equation (\ref{E:wave_data}) back in time, using the
measured data $g$ as the boundary values. The solution of this well
posed problem at $t=0$ gives the desired source function $f(x)$. Such
methods have been successfully implemented \cite{Burgh}.

Although in $2D$ or in presence of sound speed variations, Huygens'
principle does not hold anymore, and thus the signal theoretically will
stay forever, one can find good approximate solutions using a similar
approach \cite{Kuch_tat, BGK}, see discussions of the variable speed
case below.

\subsubsection{Examples of reconstructions and additional remarks about
the inversion formulas}\label{SS:examples}

\begin{itemize}

\item It is well known that different analytic inversion formulas in
tomography can behave differently in numerical implementation (e.g., in
terms of their stability), However, numerical implementation seems to
show that the analytic (backprojection type) formulas
(\ref{E:FPR3d})-(\ref{E:Kunya3D}), in spite of some of them being
not equivalent, work equally well. See, for example the results of an
analytic formula reconstruction in $3D$ shown in Fig. \ref{F:recon}.
\begin{figure}[ht]
\begin{center}
\scalebox{0.7}{\includegraphics{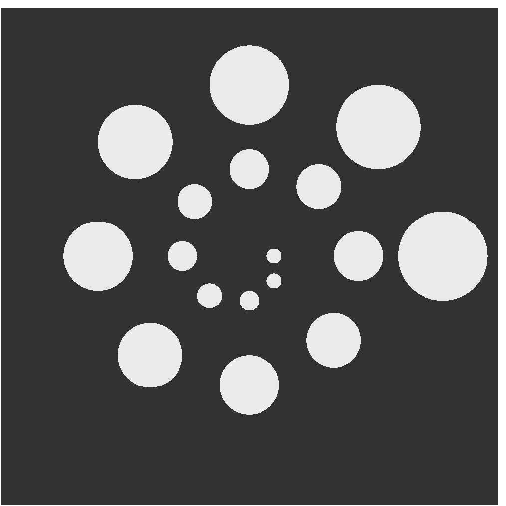}\,\,\includegraphics{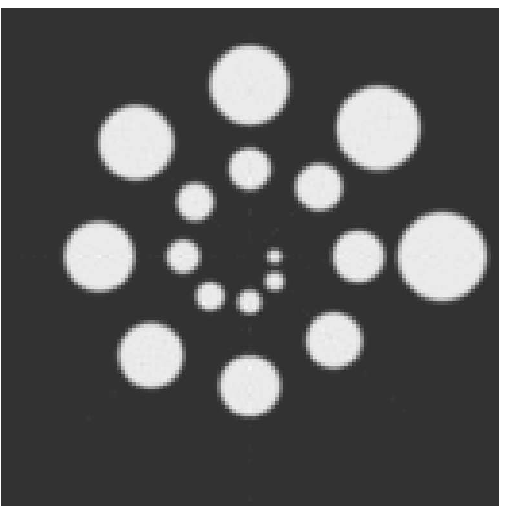}}
\end{center}
\caption{A mathematical phantom in $3D$ (left) and its reconstruction
using an analytic inversion formula.}
\label{F:recon}
\end{figure}

\item It is worth noting that although formulas
(\ref{E:FPR3d})-(\ref{E:Finch2D}) and
(\ref{E:kunyansky})-(\ref{E:Kunya3D}) will yield identical results
when applied to functions that can be represented as the spherical mean
Radon transform of a function supported inside $S$, they are in general
not equivalent when applied to functions with larger supports. Simple
examples (e.g., of $f$ being the characteristic set of a large ball
containing $S$) show that these two types of formulas provide different
reconstructions.

\item An interesting observation is that backprojection formulas
(\ref{E:FPR3d})-(\ref{E:Kunya3D}) do not reconstruct the function
$f$ correctly inside the surface $S$, if $f$ has support reaching
outside $S$. For instance, applying the reconstruction formulas to the
function $R_S (\chi_{|x|\leq 3})$ leads to an incorrect reconstruction
of the value of $f= \chi_{|x|\leq 3}$ inside $S=\{|x|\leq 1\}$. (Here
by $\chi_V$ we denote the characteristic function of the set $V$, i.e.
it takes the value $1$ in $V$ and zero outside. So, $\chi_{|x|\leq 3}$
is the characteristic function of the ball of radius $3$ centered at
the origin.)

An another example: if one adds to the phantom shown in Fig.
\ref{F:recon} two balls to the right of the surrounding sphere $S$,
this leads to strong artifacts, as seen on Fig. \ref{F:artefact}.

\begin{figure}[ht]
\begin{center}
\scalebox{0.7}{\includegraphics{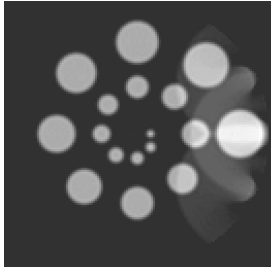}}
\end{center}
\caption{A perturbed reconstruction, due to presence of two additional
balls outside $S$ (not shown on the picture).}
\label{F:artefact}
\end{figure}

What is the reason for such a distortion? If one does not know in
advance that $f$ has support inside $S$, the backprojection formulas
shown before use insufficient information to recover a function with a
larger support, and thus uniqueness of reconstruction is lost. Then the
formulas misinterpret the data, wrongly assuming that they came form a
function supported inside $S$ and thus reconstructing the function
incorrectly.

Notice that the series reconstruction of the preceding Section is free
of such problem. E.g., the reconstruction shown in Fig. \ref{F:good}
confirms this.

\begin{figure}[ht]
\begin{center}
\scalebox{0.7}{\includegraphics{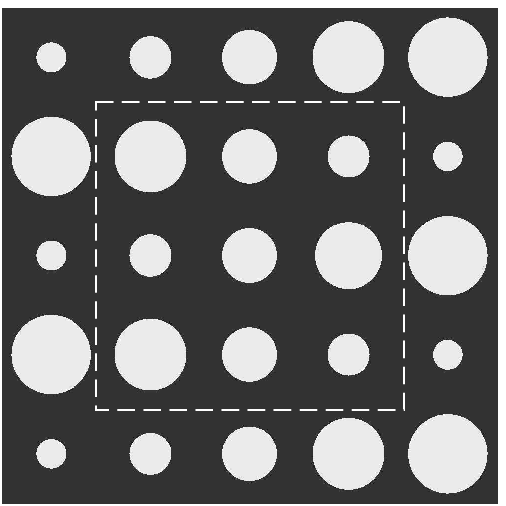}\,\,
\includegraphics{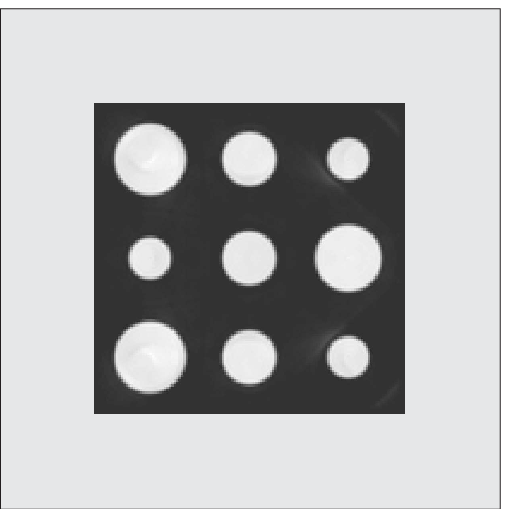}}
\end{center}
\caption{In the phantom shown on the left, most disks are
located outside the square acquisition surface $S$ indicated by the
dotted line. This does not perturb the reconstruction inside $S$
(right).}\label{F:good}
\end{figure}

\end{itemize}

\subsection{Reconstruction in the variable speed case}
\label{SS:rec_variable}
We will assume here that the sound speed $c(x)$ is smooth, positive,
constant for large $x$, and non-trapping. Although most analytic
techniques we described above do not work in the variable speed case,
some formulas can be derived and algorithms can be designed. This work
is in a beginning stage and the results described below most surely can
and will be improved.

\subsubsection{``Analytic'' inversions}\label{invers_variable}

Let us denote by $\Omega$ the interior of the observation surface $S$,
i.e. the area where the object to be imaged is located. Consider in
$\Omega$ the operator $A=-c^2(x)\Delta$ with zero Dirichlet conditions
on the boundary $S=\partial\Omega$. This operator is self-adjoint, if
considered in the weighted space $L^2(\Omega;c^{-2}(x))$.

We also denote by $E$ the operator of harmonic extension, which
transforms a function $\phi$ on $S$ to a harmonic function on $\Omega$
which coincides with $\phi$ on $S$.

The following result provides a formula for reconstructing $f$ from the
data $g$:
\begin{theorem}\cite{Kuch_tat}\label{T:reconstr_speed}
The function $f(x)$ in (\ref{E:wave_data}) can be reconstructed in
$\Omega$ as follows:
\begin{equation}\label{E:reconstruction_variable}
f(x)=(Eg|_{t=0})-\int\limits_0^\infty
A^{-\frac12} \sin{(\tau
A^{\frac12})}E(g_{tt})(x,\tau)d\tau.
\end{equation}
\end{theorem}

The validity of this result hinges upon decay estimates for the
solution (so called local energy decay \cite{Egorov, Vainb, Vainb2}),
which hold under the non-trapping condition. These estimates guarantee
a qualified decay of the solution $p(t,x)$ inside any bounded region,
e.g. in $\Omega$, when time $t$ increases. In odd dimensions decay is
exponential, but only polynomial in even dimensions. The decay can be
used instead of Huygens' principle to solve the wave equation
backwards, starting at the infinite time. This leads to the formula
(\ref{E:reconstruction_variable}).

Due to functions of the operator $A$ being involved, it is not that
clear how explicit this formula can be made. For instance, it would be
interesting to see whether one can derive from
(\ref{E:reconstruction_variable}) a backprojection inversion formula
for the case of a constant sound speed and $S$ being a sphere (we have
already seen that such formulas are known).

\subsubsection{Backpropagation}

The exponential decay at large values of time can be used as follows:
for a sufficiently large $T$, one can assume that the solution is
practically zero at $t=T$. Thus, imposing zero initial conditions at
$t=T$ and solving in reverse time direction, one arrives at $t=0$ to an
approximation of $f(x)$ \cite{BGK}.

\subsubsection{Eigenfunction expansions}

One natural way to try to use the formula
(\ref{E:reconstruction_variable}) is to use eigenfunction expansion of
the operator $A$ in $\Omega$ (assuming that such expansion is known).
This immediately leads to the following result:

\begin{theorem}\label{T:main_speed}
Under the same conditions on the sound speed as before, function $f(x)$
can be reconstructed inside
$\Omega$ from the data $g$ in (\ref{E:wave_data}), as the following
$L^2(B)$-convergent series:
\begin{equation}\label{E:coef_variable}
f(x)=\sum\limits_k f_k \psi_k(x),
\end{equation}
where the Fourier coefficients $f_k$ can be recovered using one of
the following formulas:
\begin{equation}\label{E:coef_variable2}
\begin{array}{c}
f_k=\lambda_k^{-2}g_k(0)-\lambda_k^{-3}\int\limits_0^\infty
\sin{(\lambda_k t)} g_k^{\prime\prime}(t)dt,\\
f_k=\lambda_k^{-2}g_k(0)+\lambda_k^{-2}\int\limits_0^\infty
\cos{(\lambda_k t)} g_k^{\prime}(t)dt, \mbox{ or } \\
f_k=-\lambda_k^{-1}\int\limits_0^\infty \sin{(\lambda_k t)}g_k(t)dt
=-\lambda_k^{-1}\int\limits_0^\infty \int\limits_S\sin{(\lambda_k
t)}g (x,t)\overline{\frac{\partial \psi_k}{\partial\nu}(x)}dxdt,
\end{array}
\end{equation}
and
$$
g_k(t)=\int\limits_{S}g(x,t)\overline{\frac{\partial \psi_k}
{\partial \nu}(x)}dx.
$$
Here $\nu$ denotes the external normal to $S$.
\end{theorem}

One notices that this is a generalization to the variable sound speed
case of the expansion method of \cite{Kun_series}, discussed in Section
\ref{SS:series}. An interesting feature is that, unlike in
\cite{Kun_series}, we do not need to know the whole space Green's
function for $A$ (which is certainly not known).

It is not clear yet how feasible numerical implementation of this
approach could be.

\section{Partial data. ``Visible'' and ``invisible''
singularities}\label{S:partial}

Uniqueness of reconstruction does not imply practical recoverability,
since the reconstruction procedure might be severely unstable. This is
well known to be the case, for instance, in incomplete data situations
in X-ray tomography, and even for complete data problems in some
imaging modalities, such as the electrical impedance tomography
\cite{Kuch_AMS05, KuchQuinto, Natt_old, Natt_new}.

In order to describe the results below, we need to explain the notion
of the wave front set $WF(f)$ of a function $f(x)$. This set carries
detailed information on singularities of $f(x)$. It consists of pairs
$(x,\xi)$ of a point $x$ in space and a wave vector (Fourier domain
variable) $\xi\neq 0$. It is easier to say what it means that a point
$(x_0,\xi_0)$ is {\bf not in } the wave front set $WF(f)$. This means
that one can smoothly cut-off $f$ to zero at a small distance from
$x_0$ in such a way that the Fourier transform $\widehat{\phi f}(\xi)$
of the resulting function $\phi(x)f(x)$ decays faster than any power of
$\xi$ in directions that are close to the direction of $\xi_0$. We
remind the reader that if this Fourier transform decays that way in
{\bf all} directions, then $f(x)$ is smooth near the point $x_0$. So,
the wave front set contains pairs $(x_0,\xi_0)$ such that $f$ is not
smooth near $x_0$, and $\xi_o$ indicates why it is not: the Fourier
transform does not decay well in this direction. For instance, if
$f(x)$ consists of two smooth pieces joined non-smoothly across a
smooth interface $\Sigma$, then $WF(f)$ contains pairs $(x,\xi)$ such
that $x$ is in $\Sigma$ and $\xi$ is normal to $\Sigma$ at $x$. One can
find simple introduction to the notions of microlocal analysis, such as
the wave front set, for instance in \cite{Str}.

Analysis done in \cite{Quinto} for the constant speed case
(equivalently, for the spherical mean transform $R_S$), showed which
parts of the wave front (and thus singularities) of a function $f$ can
be recovered from its partial $X$-ray data.  An analog of this result
also holds for the spherical mean transform $R_S$ \cite{LQ} (see also
\cite{XWAK} for a practical discussion). We formulate it below in an
imprecise form (see \cite{LQ} for precise formulation).

\begin{theorem}\cite{LQ}\label{T:incomplete} A wavefront set point
$(x,\xi)$ of $f$ is ``stably recoverable'' from $R_S f$ if and only if
there is a circle (sphere in higher dimensions) centered on $S$,
passing through $x$, and normal to $\xi$ at this point. \end{theorem}
As we have already mentioned, this result does not exactly hold the way
it is formulated and needs to include some precise conditions (see
\cite[Theorem 3]{LQ}). The statement is, for instance, correct if $S$
is a smooth hypersurface and the support of $f$ lies on one side of the
tangent plane to $S$ at the center of the sphere mentioned in the
theorem.

Talking about jump singularities only (i.e., interfaces between smooth
regions inside the object to be imaged), this result says that in order
for a piece of the interface to be stably recoverable (dubbed
``visible''), one should have for each point of this interface, a
sphere centered at  $S$ and tangent to the interface at this point.
Otherwise, the interface will be blurred away (even if there is a
uniqueness of reconstruction theorem). The reason is that if all
spheres of integration are transversal to the interface, the
integration smoothes off the singularity, and therefore its recovery
becomes highly unstable (numerically, one has to deal with inversion of
a matrix with exponentially fast decaying singular values). The Figure
\ref{F:incomplete} below shows an example of an incomplete data
reconstruction from spherical mean data. One sees clearly the effect of
disappearance of the parts of the boundaries that are not touched
tangentially by circles centered at transducers' locations.

\begin{figure}[ht]
\begin{center}
\scalebox{0.7}{\includegraphics{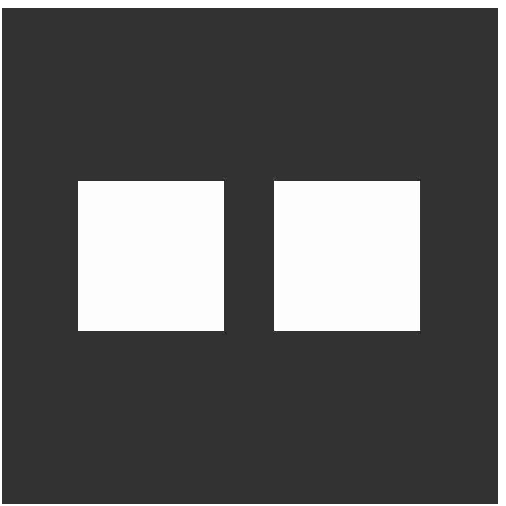}}\,
\,\scalebox{0.7}{\includegraphics{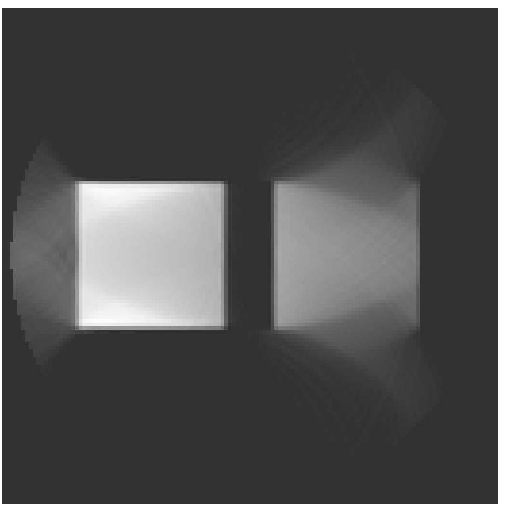}}
\end{center}
\caption{Effect of incomplete data: the phantom (left) and
its incomplete data reconstruction. The transducers were
located along a $180^o$ circular arc
(the left half of a large circle surrounding the squares).}
\label{F:incomplete}
\end{figure}

\section{Range conditions}\label{S:range}
As it has already been mentioned, the space of functions $g(t,y)$ that
could arise as exact data measured by transducers (i.e., the range of
the data), is very small (of infinite codimension in the spaces of all
functions of $t>0,y\in S$). Knowing this space (range) is useful for many
theoretical and practical purposes (reconstruction algorithms, error
corrections, incomplete data completion, etc.), and thus has attracted
a lot of attention (e.g., \cite{Leon_Radon, GGG1, GGG, GelfVil,
Helg_Radon, Helg_groups, Kuch_AMS05, KucLvin, KucLvin2, KuchQuinto, Lvin, Natt_old, Natt_new, Nessibi,
Pal_book, Q2006}.

For instance, for the standard Radon transform
$$
f(x)\to g(s,\omega)=\int\limits_{x\cdot\omega =s}f(x)dx, |\omega|=1,
$$
the range conditions on $g(s,\omega)$ are:
\begin{enumerate}
\item {\em evenness}: $g(-s,-\omega)=g(s,\omega)$

\item {\em moment conditions}: for any integer $k\geq 0$, the $k$th
moment
$$
G_k(\omega)=\int\limits_{-\infty}^{\infty}
s^k g(\omega,s)ds
$$
extends from the unit circle of vectors $\omega$ to a homogeneous
polynomial of degree $k$ in $\omega$.
\end{enumerate}
The evenness condition is obviously necessary and is kind of
``trivial''. It seems that the only non-trivial conditions are the
moment ones. However, here the standard Radon transform misleads us,
as it often happens. In fact, for more general transforms of Radon type
it is often easy (or easier) to find analogs of the moment conditions,
while analogs of the evenness conditions are often elusive (see
\cite{Kuch_AMS05, KucLvin, KucLvin2, Natt_old, Natt_new, Novikov}
devoted to the case of SPECT (single photon emission tomography)). The
same happens in TAT.

Let us deal first with the case of a constant sound speed, when one can
think of the spherical mean transform $R_S$ instead of the wave
equation model. An analog of the moment conditions was already present
implicitly (without saying that these were range conditions) in
\cite{LP1, LP2, AQ} and explicitly formulated as such in \cite{Patch}.
Indeed, our discussion in Section 4 of the polynomials $Q_k$ provides
the following conditions of the moment type:

{\bf Moment conditions} {\em \cite{AQ,LP1,LP2,Patch} on data
$g(p,r)=R_S f(p,r)$ look as follows: for any integer $k\geq 0$, the
moment
$$
M_k(\omega)=\int\limits_{0}^{\infty} r^{2k+d-1} g(p,r)dr
$$
can be
extended from $S$ to a (non-homogeneous) polynomial $Q_k(x)$ of degree
at most $2k$.}

These conditions, however, are incomplete, and in fact infinitely many
others, which play the role of an analog of evenness, need to be added.

Complete range descriptions for $R_S$ when $S$ is a circle in $2D$ were
discovered in \cite{AmbKuc_rang} and then in odd dimensions in
\cite{FR}. They were then extended to any dimension and interpreted in
several different ways in \cite{AKQ}. These conditions happen to be
intimately related to PDEs and spectral theory.

\begin{figure}[!ht]
\begin{center}
\scalebox{0.7}{\includegraphics{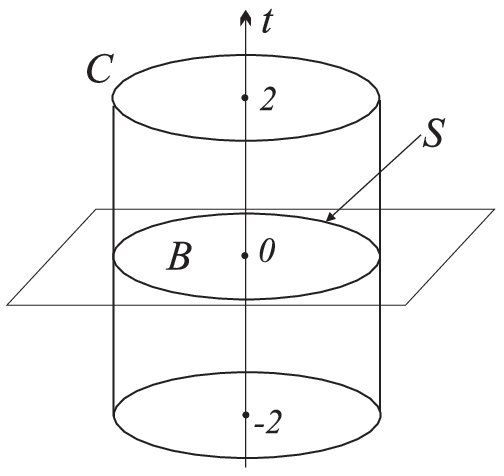}}
\end{center}
\caption{}
\label{F:}
\end{figure}

In order to describe these conditions, we need to introduce some
notations. Let $B$ be the unit ball in $\RR^d$, $S$ - the unit sphere,
and $C$ - the cylinder $B\times[0,2]$ (see Fig. \ref{F:}).

We introduce the spherical mean operator $R_S$ as before:
$$
R_Sf(x,t)=\int_{|y|=1}f(x+ty)dA(y), x\in S.
$$

Several different range descriptions for $R_S$ were provided in
\cite{AKQ}, out of which we only show a few:

\begin{theorem}\cite{AKQ}\label{T:AKQ}
    The following three statements are equivalent:
    \begin{enumerate}
    \item The function $g\in C^\infty_0 (S\times [0,2])$ is
    representable as $R_S f$ for some $f\in C^\infty_0(B)$. (In other
    words, $g$ represents an ideal (free of errors) set of TAT data.)
    \item
    \begin{enumerate}
    \item The moment conditions are satisfied.
    \item Let $-\lambda^2$ be any eigenvalue of the Laplace operator in
    $B$ with zero Dirichlet conditions and $\psi_\lambda$ be the
    corresponding eigenfunction. Then the following orthogonality
    condition is satisfied:
    \begin{equation}\label{E:orthog}
\int\limits_{S\times [0,2]} g(x,t)\partial_\nu \psi_\lambda
(x)j_{n/2-1}(\lambda t)t^{n-1}dxdt=0.
    \end{equation}
    Here $j_p(z)=c_p\frac{J_p(z)}{z^p}$ is the so called spherical
    Bessel function.
   \end{enumerate}
    \item
    \begin{enumerate}
    \item The moment conditions are satisfied.
    \item Let $\widehat{g}(x,\lambda)=\int g(x,t)j_{n/2-1}(\lambda
    t)t^{n-1}dt$. Then, for any $m\in\mathbb{Z}$, the $m^{th}$ spherical
    harmonic term $\widehat{g}_m(x,\lambda)$ of
    $\widehat{g}(x,\lambda)$ vanishes at all zeros $\lambda\neq 0$ of Bessel
    function $J_{m+n/2-1}(\lambda)$.
    \end{enumerate}
    \end{enumerate}
\end{theorem}

\begin{remark}\cite{AKQ}\label{T:AKQ2}
\begin{enumerate}
\item In odd dimensions, moment conditions are not necessary, and thus
conditions 2(b) or 3(b) suffice. (A similar earlier result was
established for a related transform in
\cite{FR}.)
\item The range conditions (2) of the previous Theorem are also
necessary when $S$ is the boundary of any bounded domain, not
necessarily a sphere.
\item An analog of these conditions can be derived for a variable sound
speed (without non-trapping conditions imposed).
\end{enumerate}
\end{remark}

\section{Concluding remarks}\label{S:remarks}

\subsection{Variations of the TAT procedure}

\subsubsection{Planar and linear transducers}\label{SS:linear and planar}

Assuming that transducers are point-like, is clearly an approximation,
and in fact a transducer measures the average pressure over its
area. It has been rightfully claimed that the point approximation for
transducers should lead to some blurring in the reconstructions. This,
as well as intricacies of reconstructions from the data obtained by
point transducers, triggered recent proposals for different types of
transducers (see \cite{haltmaier_int_and line, haltmaier_fabri},
\cite{haltmaier_large}-\cite{haltmaier_challenge_line},
\cite{haltmaier_opt_line, haltmaier_interf}). In these papers, it was
suggested to use either planar, or line detectors.

In the first case \cite{haltmaier_large}, the detectors are assumed to
be large and planar, ideally assumed to be approximations of infinite
planes that are placed tangentially to a sphere containing the object.
Thus, the data one collects is the integrals of the pressure over these
planes, for all values of $t>0$. If one takes the standard $3D$ Radon
transform of the pressure $p(x,t)$ with respect to $x$:
$$
P(x,t)\mapsto q(s,t,\omega)=\int\limits_{x\cdot\omega=s}p(x,t)dA(x),
$$
where $dA$ is the surface measure and $\omega$ is a unit vector in
$\RR^3$, this is well known to reduce the $3D$ Laplace operator
$\Delta_x$ to the second derivative $\partial^2/\partial s^2$
\cite{Leon_Radon, GGG1, GGG, GelfVil,Helg_Radon, Helg_groups}, and thus
the $3D$ wave equation to the string vibration problem. The measured
data provide the boundary conditions for this problem. The  initial
conditions in (\ref{E:wave}) mean evenness with respect to time, and
thus the standard d'Alambert formula leads to the immediate realization
that the measured data is just the $3D$ Radon transform of $f(x)$.
Thus, the reconstruction boils down to the well known inversion
formulas for the Radon transform.

Another proposal (\cite{haltmaier_int_and line, haltmaier_fabri},
\cite{haltmaier_int_line}-\cite{haltmaier_challenge_line},
\cite{haltmaier_opt_line, haltmaier_interf}) is to use line detectors
that provide line integrals of the pressure $p(x,t)$. Such detectors
can be implemented optically, using either Fabry-Perot
\cite{haltmaier_int_and line}, or Mach-Zehnder \cite{haltmaier_interf}
interferometers.

Suppose that the object is surrounded by a surface that is rotation
invariant with respect to the $z$-axis. It is suggested to place the
line detectors perpendicular to the $z$-axis and tangential to the
surface. The same consideration as above then shows that after the $2D$
Radon (or $X$-ray, which in $2D$ is the same) transform in each plane
orthogonal to $z$-axis, the $3D$ wave equation converts into the $2D$
one for the Radon data. The measurements provide the boundary data.
Thus, the reconstruction boils down to solving a $2D$ problem similar
to the one in the case of point detectors, and then inverting the $2D$
Radon transform.

Due to the recent nature of these two projects, it appears to be too
early to judge which one will be superior in the end. For instance, it
is not clear beforehand, whether the approximation of infinite size
(length, area) of the linear or planar detectors works better than the
zero dimension approximation for point detectors. Further developments
will resolve these questions.

\subsubsection{Direct imaging techniques}

Some direct imaging techniques have been suggested, which might not require
mathematical reconstructions. See, for instance, \cite{Freehand} about an
acoustic lens system.

\subsubsection{Using contrast agents}

Contrast agents to improve TAT imaging have been developed (e.g.,
\cite{Copland}).

\subsubsection{Passive thermoacoustic imaging}

The TAT model we have considered can be called ``active thermoacoustic
tomography,'' due to the set-up when the practitioner creates the
signal. There has been some recent development of the ``passive
thermoacoustic tomography,'' where the thermoacoustic signal is used to
image the temperature sources present inside the body. One can find a
survey of this area in \cite{passive}.

\subsection{Uniqueness}\label{SS:remunique}

\subsubsection{Sketch of the proof of Theorem \ref{T:AKQ}}
We provide here a brief outline of the rather technical proof of
Theorem \ref{T:AKQ}.

Suppose that $f$ is compactly supported, not identically zero, and such
that $R_S f=0$. Our previous considerations show that one can assume
that $S$ is an algebraic curve (not a straight line) that is contained
in the set of zeros of a non-trivial harmonic polynomial. Now one
touches the boundary of the support of $f$ from outside by a circle
centered on $S$. Then microlocal analysis of the operator $R_S$ (which
happens to be an analytic Fourier Integral Operator, FIO \cite{Be,
GU91, Guill75, Guill85, GS, KLM, Q1980}) shows that, due to the equality
$R_Sf=0$, at the tangency point the vector co-normal to the sphere
should not belong to the analytic wave front of $f$ (microlocal
regularity of solutions of $R_S f=0$). This, for instance, can be also
extracted from the results of \cite{StefUhl}. On the other hand, a
theorem by H\"{o}rmander and Kashiwara \cite[Theorem 8.5.6]{Horm} shows
that this vector must be in the analytic wave front set, since $f=0$ on
one side of the sphere (a microlocal version of uniqueness of analytic
continuation). This way, one gets a contradiction. Unfortunately, the
life is not so easy, and the proof sketched above does not go through
smoothly, due to possible cancelation of wavefronts at different
tangency points. Then one has to involve the geometry of zeros of
harmonic polynomials \cite{Flatto} to exclude the possibility of such a
cancelation.

Thus, the proof uses microlocal analysis and geometry of zeros of
harmonic polynomials. Both these tools have their limitations. For
instance, the microlocal approach (at least, in the form it is used in
\cite{AQ}) does not allow considerations of non-compactly supported
functions. Thus, the validity of the Theorem for arbitrarily fast
decaying, but not compactly supported, functions is still not
established, albeit it most certainly holds. On the other hand, the
geometric part does not work that well in dimensions higher than two.
Development of new approaches is apparently needed in order to overcome
these hurdles. A much simpler PDE approach has emerged recently
\cite{FPR} (see also \cite{AmbKuc_inj} and the next Section), albeit
its achievements have been limited so far.

\subsubsection{Some open problems concerning uniqueness}

As it has already been mentioned, one can consider the practical
problems about uniqueness resolved. However, the mathematical
understanding of the uniqueness problem for the restricted spherical
mean operators $R_S$ is still unsatisfactory. Here are some questions
that still await their resolution:

\begin{enumerate}

\item Describe uniqueness sets in dimensions larger than $2$ (prove the
Conjecture \ref{C:n-dim}). Recent limited progress, as well as
variations on this theme can be found in \cite{Agr}-\cite{AmbKuc_inj}.

\item Prove Theorem \ref{T:AQ} without using microlocal and harmonic
polynomial tools.

\item Prove Theorem \ref{T:AQ} on uniqueness sets $S$ under
the condition of sufficiently fast decay (rather than
compactness of support) of the function. Very little
is known for the case of functions without compact
support. The main known result is of \cite{ABK}, which
describes for which values of $1\leq p \leq \infty$ the
result of Corollary \ref{C:closed} still holds:

\begin{theorem}\cite{ABK}
Let $S$ be the boundary of a bounded domain in $\RR^d$ and $f\in
L^p(\RR^d)$ such that $R_S f\equiv 0$. If $p\leq 2d/(d-1)$, then
$f\equiv 0$ (and thus $S$ is injectivity set for this space). This
fails for any $p>2d/(d-1)$.
\end{theorem}

\end{enumerate}

\subsection{Inversion}\label{SS:reminversion}
Albeit closed form (backprojection type) inversion formulas are
available now for the cases of $S$ being a plane (and object on one
side from it), cylinder, and a sphere, there is still some mystery
surrounding this issue.

\begin{enumerate}
\item Can one write a backprojection type inversion formula in the case
of the constant sound speed for a closed surface $S$ which is not a
sphere? We suspect that the answer to this question is negative (see also
related discussion in \cite{Anastasio, Beylkin}).

\item The inversion formulas for $S$ being a sphere assume that the
object to be imaged is inside $S$. One can check on simplest examples
that if the support of function $f(x)$ reaches outside $S$, the
inversion formulas do not reconstruct the function correctly even
inside of $S$. See \cite{AKK} for a discussion.

\item The I.~Gelfand's school of integral geometry has developed a
marvelous machinery of the so called $\kappa$ operator, which provides
a general approach to inversion and range descriptions for transforms
of Radon type \cite{GGG1, GGG}. In particular, it has been applied to
the case of integration of various collections (``complexes'') of
spheres in \cite{GGG, Gi}. This consideration seems to suggest that one
should not expect explicit closed form inversion formulas for $R_S$
when $S$ is a sphere. We, however, know that such formulas have been
discovered recently \cite{FPR, Kunyansky}. This apparent controversy
has not been resolved.

\item Can one derive any more explicit analytic formulas from
(\ref{E:reconstruction_variable})?

\item Can the series expansion formulas of Theorem \ref{T:main_speed}
be efficiently implemented?

\end{enumerate}

One can also mention that in some works \cite{Anastasio, Klibanov} it is
suggested to use in the TAT problem not only the values of the pressure
measured by transducers on the observation surface $S$, but its normal
derivative to $S$ as well. If one knows both, then taking Fourier transform
in the time variable and using the whole space Green's function for the
Helmholtz equation leads immediately to a reconstruction formula for the
solution (which seems to be much simpler than what is proposed in
\cite{Klibanov}). The problem is that this normal derivative is not measured
by TAT devices. Under some circumstances (e.g., when there are no sources of
ultrasound outside $S$), one can prove the theoretical possibility of
recovering the missing normal derivative. This, however, does not seem to us
to be a plausible procedure. In rare cases (planar, cylindrical, or
spherical surface $S$), when involvement of the normal derivative can be
eliminated (e.g., \cite{Anastasio, Beylkin}), this might lead to feasible
inversion algorithms, but in these cases, as explained before in this text,
explicit and nicely implementable analytic inversion formulas are available.
So, jury is still out on this issue as well.

\subsection{Stability}\label{SS:remstability}

Stability of inversion when $S$ is a sphere surrounding the support of
$f(x)$ is the same as for the standard Radon transform, as the results
of \cite{Palam_funk} and second statement of Theorem \ref{T:AKQ2} show.
However, if the support reaches outside, albeit Corollary
\ref{C:closed} still guarantees uniqueness of reconstruction, stability
(at least for the parts outside $S$) is gone. Indeed, Theorem
\ref{T:incomplete} shows that some parts of singularities of $f$
outside $S$ will not be stably ``visible.''

\subsection{Range}\label{SS:remrange}

As Theorem \ref{T:incomplete} states, the range conditions 2 and 3 of
Theorem \ref{T:AKQ} are necessary also for non-spherical closed
surfaces $S$ and for functions with support outside $S$. They, however,
are not expected to be sufficient, since  Theorem \ref{T:incomplete}
indicates that one might expect non-closed ranges in some cases. The
same applies for non-constant sound speed case.

\section*{Acknowledgments}\label{S:acknowledgments}

The work of the first author was partially supported by the NSF DMS
grants 0604778 and 0648786. The second author was partially supported
by the DOE grant DE-FG02-03ER25577 and NSF DMS grant 0312292. Part of
this work was completed when the first author was at the Isaac Newton
Institute for Mathematical Sciences. The authors express their
gratitude to the NSF, DOE and INI for this support. The authors thank
M.~Agranovsky for extremely useful discussions, M.~Anastasio, G.~Beylkin
and M.~Klibanov for providing preprints and references, and the reviewers for
very helpful remarks.


\end{document}